\documentclass[12pt]{article}
\usepackage{amsxtra,amssymb,amsthm,amsmath,latexsym}

\theoremstyle{plain}

\def\R{{\mathbb R}}

\def\oH{{\overset{\circ}{H}}}
\def\oH1{{\overset{\circ}{H}\kern-.02in{}^1}}

\def\bee{\begin{equation*}}
\def\eee{\end{equation*}}
\def\be{\begin{equation}}
\def\ee{\end{equation}}

\begin{document}

\title{Completeness of the set $\{e^{ik\beta \cdot s}\}|_{\forall \beta \in S^2}$
}

\author{Alexander G. Ramm\\
 Department  of Mathematics, Kansas State University, \\
 Manhattan, KS 66506, USA\\
ramm@math.ksu.edu\\
http://www.math.ksu.edu/\,$\sim $\,ramm
}

\date{}
\maketitle\thispagestyle{empty}

\begin{abstract}

It is proved that the set   $\{e^{ik\beta \cdot s}\}|_{\forall \beta \in S^2}$, where $S^2$ is the unit sphere in $\mathbb{R}^3$,  $k>0$ is a fixed
constant, $s\in S$, 
is total in $L^2(S)$ if and only if $k^2$ is not a Dirichlet eigenvalue of the Laplacian in $D$.  Here $S$ is a smooth, closed, connected surface in  $\mathbb{R}^3$.
\end{abstract}

\section{Introduction}\label{S:1}
 Let $D\subset \mathbb{R}^3$ be a bounded domain with a connected closed $C^2-$smooth boundary $S$,
 $D':=\mathbb{R}^3\setminus D$ be the unbounded exterior domain and $S^2$ be the unit sphere in $\mathbb{R}^3$, $\beta \in S^2$, $s\in S$.

We are interested in the following problem:

{\em Is the set $\{e^{ik\beta \cdot s}\}|_{\forall \beta \in S^2}$ total in $L^2(S)$?}

A set $\{\phi(s,\beta)\}$ is total (complete) in $L^2(S)$ if the relation $\int_S f(s)\phi (s,\beta)ds=0$ for all $\beta \in S^2$ implies $f=0$,
where $f\in L^@(S)$ is an arbitrary fixed function.

The above question is of interest by itself, but also it is of interest in scattering problems and in inverse problems, see \cite{R190}--\cite{R661}.

Our result is: 

{\bf Theorem 1.}  {\em The set  $\{e^{ik\beta \cdot s}\}|_{\forall \beta \in S^2}$ is total in $L^2(S)$
if and only if $k^2$ is not a Dirichlet eigenvalue of the Laplacian in $D$.}
\vspace{3mm}

\section{Proof of Theorem 1}\label{S:2}

{\bf Necessity.}  Let $f\in L^2(S)$ and 
\be\label{e1}
\int_S f(s)e^{ik\beta \cdot s}ds=0\quad \forall \beta \in S^2,
\ee
and there is a $u\not\equiv 0$ such that 
\be\label{e2}
(\nabla^2+k^2)u=0 \quad in\quad D, \qquad u|_{S}=0.
\ee
Choose $f=u_N$, where $N$ is the unit normal to $S$ pointing out of $D$. Then, by Green's formula,
equation \eqref{e1} holds and $f\not\equiv 0$
by the uniqueness of the solution to the Cauchy problem for elliptic equation \eqref{e2}. Necessity
is proved.

{\bf Sufficiency.} Assume that $f\in L^2(S)$ is
and arbitrary fixed function, $f\not\equiv 0$,
and \eqref{e1} holds. Let $h\in L^2(S^2)$ be arbitrary and 
\be\label{e3}
w(x):=\int_{S^2}h(\beta)e^{ik\beta \cdot x}.
\ee
Then 
\be\label{e4} 
(\nabla^2+k^2)w=0 \quad in\quad \R^3.
\ee
If \eqref{e1} holds, then
\be\label{e5} 
\int_S f(s)w(s)ds=0
\ee
for all $w$ of the form \eqref{e3}.
Let us now apply the following Lemma:

{\bf Lemma 1.} {\em The set $\{w|_S\}$ for all $h\in L^2(S^2)$ is the orthogonal complement in $L^2(S)$ to the linear span of the set $\{v_N\}$,
where $v$ solve equation \eqref{e4} and $v|_S=0$.}

If $k^2$ is not a Dirichlet eigenvalue of the Laplacian in $D$, then Lemma 1 implies that the set 
 $\{w|_S\}$ is total in $L^2(S)$, so  \eqref{e1}
 implies $f=0$. Sufficiency and Theorem 1 are proved. \hfill$\Box$

Lemma 1 is similar to Theorem 6 in \cite{R666}. 

{\bf Proof of Lemma 1.}  Let $w|_S:=\psi$. Choose an arbitrary $F\in C^2(D)$ such that $F|_S=\psi$.
Define $G:=F-w$ in $D$. Then
\be\label{e6}
(\nabla^2+k^2)G=(\nabla^2+k^2)F \quad in D;\quad G|_S=0.
\ee
For \eqref{e6} to hold it is necessary and sufficient that
\be\label{e7}
0=\int_D (\nabla^2+k^2)F v dx,
\ee
where $v$ is an arbitrary function in the set
of solutions of equation \eqref{e2}. Using
Green's formula one reduces condition \eqref{e7}
to the following condition:
\be\label{e8}
\int_S \psi v_N ds=0.
\ee
Therefore the set $\{\psi\}$ is the orthogonal
complement in $L^2(S)$ of the linear span of the functions $\{v_N\}$. Lemma 1 is proved. \hfill$\Box$


\begin{thebibliography}{1000} 




\bibitem{R190} A.G.Ramm,  {\em Scattering by obstacles}, D.Reidel, Dordrecht, 1986.

\bibitem{R470} A.G.Ramm, {\em Inverse problems}, Springer, New York, 2005.

\bibitem{R666} A.G.Ramm,  Solution to the Pompeiu problem and the related symmetry problem,

Appl. Math. Lett., 63, (2017), 28-33.

\bibitem{R672} A.G.Ramm, Perturbation of zero surfaces, {\em Global Journ. of Math. Analysis},
5, (1),  (2017), 27-28. 

\bibitem{R661} A.G.Ramm, Uniqueness of the solution to inverse obstacle scattering with
non-over-determined data, Appl. Math. Lett., 58, (2016), 81-86.

\end{thebibliography}
\end{document}